\documentclass[lettersize,journal]{IEEEtran}
\usepackage{amsmath,amsfonts,amsthm}
\usepackage{algorithmic}
\usepackage{algorithm}
\usepackage{array}
\usepackage[caption=false,font=normalsize,labelfont=sf,textfont=sf]{subfig}
\usepackage{textcomp}
\usepackage{stfloats}
\usepackage{url}
\usepackage{siunitx}
\usepackage{verbatim}
\usepackage{graphicx}
\usepackage{cite}
\usepackage{amssymb, dsfont}
\usepackage{cleveref}
\usepackage{xparse}
\usepackage{makecell}
\newtheorem{definition}{\textbf{Definition}}

\newtheorem{assumption}{\textbf{Assumption}}

\begin{document}
    \title{Absolute Ranking: An Essential Normalization for Benchmarking Optimization Algorithms}

    \author{
        Yunpeng Jing, Qunfeng Liu,
        \thanks{Yunpeng Jing, Qunfeng Liu are with the School of Computer Science and Technology,
        Dongguan University of Technology, Dongguan 523808, China.}
    }

    \markboth{IEEE Transactions on Evolutionary Computation}
    {Shell \MakeLowercase{\textit{et al.}}: Absolute Ranking: An Essential Normalization for Benchmarking Optimization Algorithms}


    \maketitle

    \begin{abstract}
        Evaluating performance across optimization algorithms on many problems presents a complex challenge due to the diversity of numerical scales involved.
        Traditional data processing methods, such as hypothesis testing and Bayesian inference, often employ ranking-based methods to normalize performance values across these varying scales.
        However, a significant issue emerges with this ranking-based approach: the introduction of new algorithms can potentially disrupt the original rankings.
        This paper extensively explores the problem, making a compelling case to underscore the issue and conducting a thorough analysis of its root causes.
        These efforts pave the way for a comprehensive examination of potential solutions.
        Building on this research, this paper introduces a new mathematical model called ``absolute ranking'' and a sampling-based computational method.
        These contributions come with practical implementation recommendations, aimed at providing a more robust framework for addressing the challenge of numerical scale variation in the assessment of performance across multiple algorithms and problems.

    \end{abstract}


    \section{Introduction}\label{sec:introduction}
    \IEEEPARstart{A}{lgorithm} benchmarking is of paramount importance for the development of optimization algorithms~\cite{carrasco2020recent,halim_performance_2021,abdel2018metaheuristic}, as it provides the foundation for algorithm selection and parameter optimization.
The primary means of conducting algorithm benchmarking currently involve applying algorithms to benchmark test sets, conducting experiments, and analyzing the results.
However, analyzing experimental results faces a challenge: different benchmark functions have different numerical scales.
A common approach to address the issue is to normalize the data generated under each benchmark function.
For instance, non-parametric hypothesis testing methods~\cite{garcia_study_2009} employ ranks for normalization, while Bayesian inference~\cite{benavoli2017time,rojas-delgado_bayesian_2022}, considered as an alternative, uses pairwise comparisons to mitigate the influence of problem scale.

However, this paper has identified a potential issue with the use of rank-based or pairwise comparison-based methods.
Specifically, an algorithm's rank or its odds of winning pairwise comparisons may be influenced by the composition of the algorithm set.
This means that the comparison result of two algorithms may change due to the addition or removal of irrelevant third-party algorithms, contradicting the intuitive expectation that comparison result of algorithms should depend solely on each algorithm's performance on the test set and be independent of other irrelevant algorithms.
This phenomenon is referred to in this paper as Non-Independence from Irrelevant Alternatives, abbreviated as NIIA\@.

To validate the existence of this phenomenon in algorithm benchmarking, the paper designs two datasets to test whether an analysis method exhibit NIIA problems.
The first dataset is designed as a subset of the second dataset, and both of them are analyzed by the same method.
If the results from the second dataset contradict compared to the results from the first dataset, it indicates the method exhibit NIIA problem.
The experimental results reveals that both non-parametric hypothesis testing and Bayesian inference exhibit NIIA problem.

In response to the NIIA problem, this paper introduces a new concept called ``absolute rank''.
It can be understood as a scheme that includes all possible algorithms for data analysis in the algorithm set.
By this way, the presence or absence of finite number of algorithms does not affect the value of the absolute rank, thus avoiding the NIIA problem.
This paper provides a mathematical model of the absolute ranking and discusses both analytical methods for calculating an exact solution and using sampling methods to estimate an approximate solution.

In addition to theoretical work, this paper also discusses practical applications, offering recommendations on how to address the NIIA problem under low-cost conditions and presents a simple case to illustrate the process.

The subsequent chapters of this paper are organized as follows: Section~\ref{sec:preliminaries} introduces some preliminaries and clarifies the motivation for this paper in more detail, Section~\ref{sec:the-issue-of-relative-ranking} expounds upon the existence of the NIIA issue and scrutinizes its underpinnings, Section~\ref{sec:absolute-ranking} delves deeply into the pertinent mathematical models, definitions, and calculations of the proposed solution, Section~\ref{sec:applications} provides a discussion of real-world applications of the absolute ranking, and finally, Section~\ref{sec:conclusions} encapsulates the conclusion and provides a forward-looking perspective.

    \section{Preliminaries}\label{sec:preliminaries}
    In this section, we progressively introduce the background knowledge of algorithm benchmarking, present various symbols and conventions, and ultimately lead to the main focus of this paper.

\subsection{Conventions for Notation}\label{subsec:conventions-for-notation}
To prevent excessive verbosity in descriptions, avoid the introduction of too many concepts and redundant definitions, default symbol conventions will be adopted.
Most of these conventions typically do not lead to misunderstandings, however, for the sake of rigor, some of the rules followed in this paper are explicitly listed below.

\subsubsection{partial function}
If there is a function $g(x,y,z)$ that is a function with three independent variables, then $g_x^y(z)=g(x,y,z)$ is a function with a single variable $z$, where $x$ and $y$ are fixed.
When the values of $x$ or $y$ can be inferred from the context, they may be omitted, and the determination of which parameters are omitted is based on the data types of the arguments.
The standalone function name is sometimes used to represent the return value of the function, like $g=g(x)$.

\subsubsection{composite function}
If there is a function $g(x)$ whose independent variables can be derived from formula $x=h(u,v)$, then $g(u,v)$ can be used to represent $g(h(u,v))$, if $h$ has been provided in the context.

\subsubsection{parallel mapping}
Let $X$ be a matrix (or vector, set, high-dimensional array, random variable) with all its elements within the domain of function $g$, then $g(X)$ returns data that is structurally identical to $X$, but with the respective elements mapped to their corresponding function values.
The word ``set'' in this paper is assumed to represent a "multiset".
This is done to ensure that when a set $X$ is mapped, its cardinality (number of elements) remains unchanged, i.e., $|g(X)| = |X|$.

\subsubsection{lambda expression}
Lambda expressions can be used to represent an anonymous function.
When a function's return value is not a specific value but an anonymous function, we can use a lambda expression to represent the return value of this function.
For example, let $g(x)=(\lambda)\to\lambda+x$, that means the function $g$ takes $x$ as its parameter and returns another function, which takes $\lambda$ as its parameter and returns $\lambda + x$.
In this case, $g(2)$ is a function that can be called, and $g(2)(3)$ returns 5.

\subsection{The search sequences of the algorithms}\label{subsec:the-search-sequences-of-the-algorithms}
This paper addresses the typical scenario of using metaheuristic algorithms to solve black-box single-objective optimization problems.
A fundamental concept in the this field is famous ``No Free Lunch Theorem for Optimization'', which states that no single optimization algorithm or strategy is universally superior for all types of optimization problems.
The practical implication is that algorithm benchmarking should be problems-specific.
This also implies that a benchmark function set used for algorithm benchmarking, although claimed to be black-box, inherently possesses certain characteristics, otherwise it cannot be used to evaluate the performance of algorithms.
The effectiveness of algorithms, if they outperform random search, suggests that their search behavior to some extent aligns with these characteristics.

\begin{figure}[!t]
    \centering
    \includegraphics[width=3.5in]{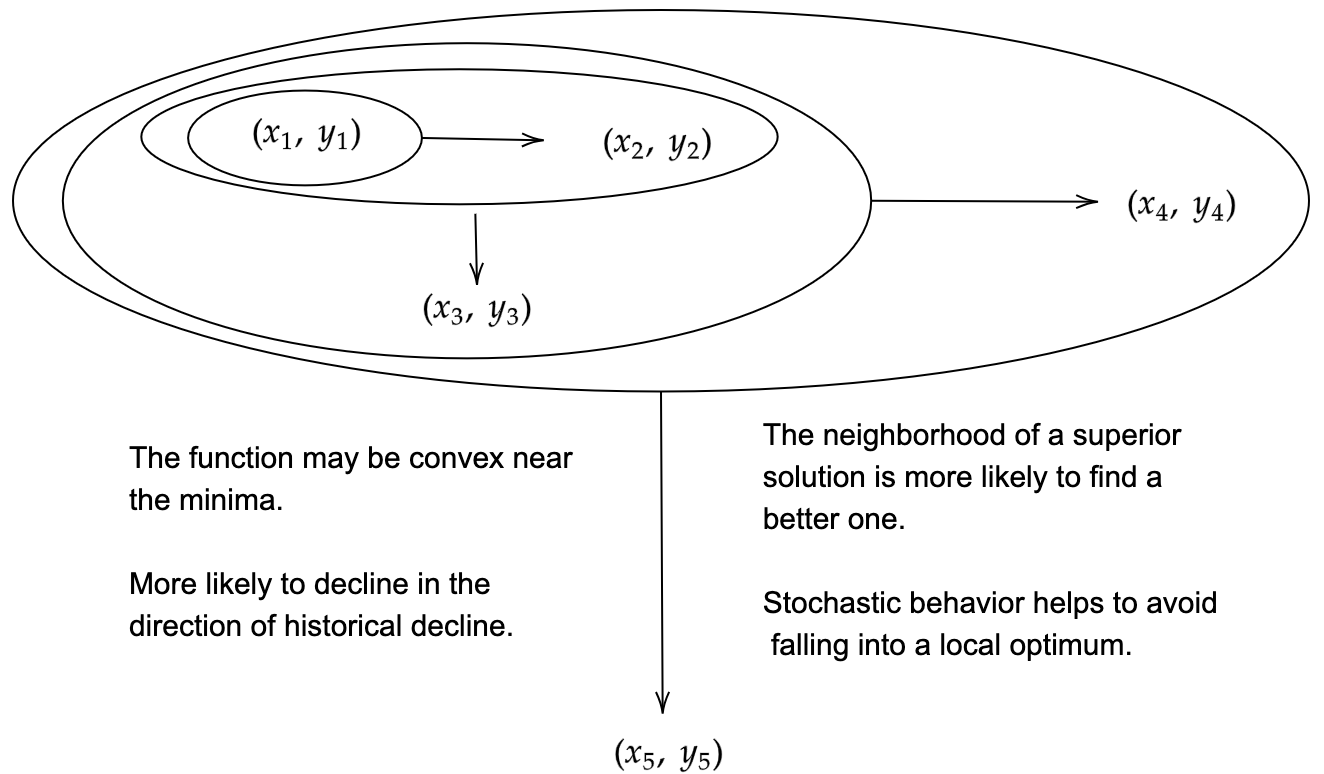}
    \caption{Typical search behavior of metaheuristic algorithms.
    The next search target depends on the historical search results and the algorithm's meta-heuristic inference strategy.}
    \label{fig:search}
\end{figure}

Metaheuristic algorithms typically do not explicitly specify the characteristics of the optimization problems they are applicable to.
However, they rely on heuristic inferences that allow them to deduce which points are more likely to yield smaller values, with the help of historical search results.
This inference guides their next search direction, leading to the next round of function evaluations iteratively, until their computational resources are exhausted.
The Fig.~\ref{fig:search} displays an illustration of the search process.

The term ``computational resources'' here refers to the maximum allowed number of calls to the objective function that an algorithm can make, the maximum number is denoted as $c$.
In this context, the search sequence of a single run of algorithm $a$ on a benchmark function $f$ can be represented as follows:
\begin{equation}
    s = s(a, f) = (x_1, x_2, \dots, x_n)
\end{equation}
The corresponding values of the objective function that are found during the search are represented as:
\begin{equation}
    f(x) = (f(x_1), f(x_2), \dots, f(x_n)).
\end{equation}

In practical experimental scenarios, there may be some minor differences in details.
For example, when an algorithm prematurely determines the optimal value of the function and reaches a predetermined level of accuracy before computational resources are exhausted.
Or, the algorithm iterates in the form of a population.
This paper overlooks these minor differences and treats them as step-by-step calls until computational resources are depleted.

\subsection{Performance metrics of the algorithms}\label{subsec:performance-metrics-of-the-algorithms}
The word ``metric'' in algorithm benchmarking can have different meanings in different contexts.
In this paper, it refers to a numerical evaluation of an algorithm on a specific problem.
For deterministic algorithms, the performance metric is typically a function of $f(s)$,
\begin{equation}
    m^* = m^*(f(s)) = m^*(s).
\end{equation}

In order to avoid getting stuck in local optima or preventing potential deadlocks, algorithms introduce some random behavior.
If the random functions used by the algorithm are deterministic, we assume that the search sequences also follow a deterministic distribution, and consequently, the corresponding performance metrics also follow a deterministic distribution.
A statistic of this distribution, such as the mean or median, can be used to assess the performance of non-deterministic algorithms.
In practice, we often estimate this statistic by conducting multiple repeated experiments.
If the experiments are repeated $r$ times, each round resulting in a search sequence $s$, then we obtain $r$ sequences, form a search matrix $X_{r*c}$.
The corresponding function values of $X$ is $f(X)$, the two dimensions of which are referred to as the cost-dimension and the round-dimension.
The performance metric $m$ for non-deterministic algorithms can be considered a function that takes $f(X)$ as input and returns a statistic value.
\begin{equation}
    m = m(f(X)) = m(X).
\end{equation}
The performance metric involves two successive dimensionality reduction operations on the original two-dimensional data $f(X)$.
Firstly, in the cost-dimension, each $f(s)$ is reduced to its corresponding $m^*$.
Then, in the round-dimension, $r$ values of $m^*$ are further reduced to $m$ through a statistical measure.

To avoid redundant classification discussions, deterministic algorithms will be considered as the case where $r=1$, or the case where $r$ experiments yield identical results.
Consequently, this paper will no longer distinguish between deterministic and non-deterministic algorithms in the following sections.

Assuming there are $n$ algorithms, denoted as $\mathbb{A}=\{a_1, a_2, \dots, a_n\}$, and $p$ benchmark functions, denoted as $\mathbb{F}=\{f_1, f_2, \dots, f_p\}$.
After conducting performance metric analysis for each algorithm-function pair, we obtain a two-dimensional performance matrix $M_{n*p}=m(\mathbb{A}*\mathbb{F})$, where $M_{ij}$ is the performance metric for algorithm $a_i$ on function $f_j$.
From this matrix, we can easily compare the performance metrics of two algorithms on the same benchmark function.
Furthermore, when we want to evaluate an algorithm's performance across the problem set $\mathbb{F}$, the rows in the $M$ (performance across different benchmark functions for the same algorithm) should be aggregated.
However, the differing benchmark problems can impact the numerical scale of performance metrics, posing challenges for aggregation across the problem dimension.

\subsection{Different scales: the challenge in aggregating across the problem-dimension.}\label{subsec:different-scales:-the-challenge-in-aggregating-across-the-problem-dimension.}
\begin{figure}[!t]
    \centering
    \includegraphics[width=3.5in]{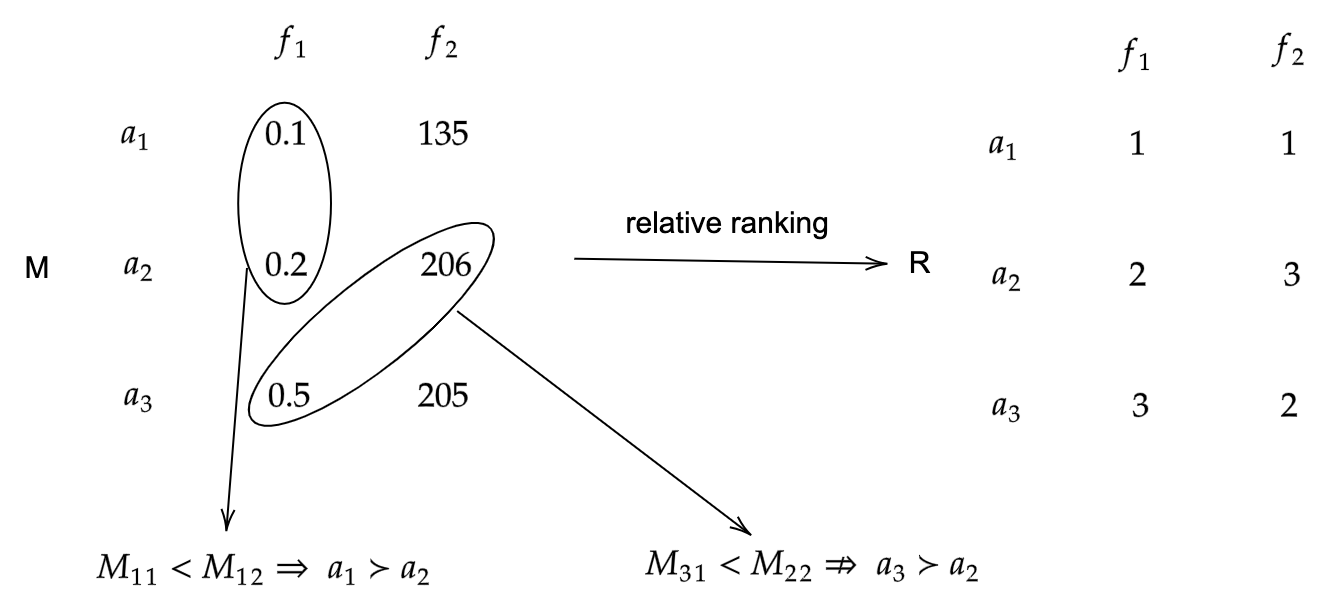}
    \caption{Illustration of different scales and rank normalization.
    When the metrics for different problems exhibit differences in numerical scales, the numerical values across algorithms become non-comparable.
    Rank normalization is one way to address this issue, as it transforms the numerical values for all problems to the range from 1 to $n$.}
    \label{fig:challenge}
\end{figure}

The matrix on the left in Fig.\ref{fig:challenge} provides us with an intuitive example of different data scales.
It illustrates that the characteristics of the benchmark problems themselves may also influence algorithm performance metrics.
This influence can be negligible when multiple algorithms are solving the same problem, adhering to the principle of controlling variables.
However, when analyzing numerical values across different problems, this influence cannot be disregarded.
Due to the varying data scales, not only can these data points not be compared numerically, but their statistical measures, such as mean, median, variance, and so forth, also lose their meaning.

To make data aggregation meaningful, we should transform the data to a common scale.
Let $M_{*j}$ be the $j$-th column of $M$, the two common normalizations in statistics, Max-Min Scaling and z-score, are represented as,
\begin{equation}
    \label{eq:mms}
    MMS(M_{ij}) = \frac{M_{ij}-\min\{M_{*j}\}}{\max\{M_{*j}\}-\min\{M_{*j}\}},
\end{equation}
and
\begin{equation}
    \label{eq:z-score}
    z(M_{ij}) = \frac{M_{ij}-\mu}{\sigma},
\end{equation}
where $\mu$ is the mean value of $M_{ij}$ and $\sigma$ is the standard deviation of $M_{*j}$.

The commonly used data analysis methods, non-parametric hypothesis testing, and Bayesian inference, all employ (or indirectly use) rank normalization to address this issue, which will be discussed in the following subsections.

\subsection{Rank normalization and relative ranking}\label{subsec:rank-normalization-and-relative-ranking}
As shown on the right part of Fig.~\ref{fig:challenge}, rank normalization transforms the elements in $M$ into their order within their respective columns.
Let $\pi_S(e)$ represent the order of $e$ within set $S$ (with the smallest assigned a value of 1 and the largest assigned a value of $|S|$).
The result of rank normalization of $M$, which is denoted as $R_{n*p}$, is calculated by,
\begin{equation}
    R_{ij} = \pi_{M_{*j}}(M_{ij})
\end{equation}
After rank normalization, each column of the data $R$ contains values ranging from 1 to $n$.
These values are no longer influenced by the characteristics of the benchmark functions in terms of scale but are determined solely by the performance of the algorithms.

The results of rank normalization may have various equivalent variations, such as simultaneously dividing all values by $n$ to ensure that all values are within $[0, 1]$, or splitting the complete sorting of a column into ${n\choose 2}$ pairwise comparison results, as Bayesian inference done.
These variants, although different in form, essentially maintain consistency.
Based on this, this paper introduces the definition of relative ranking:
\begin{definition}[relative ranking]
    A data analysis method is considered to utilize relative ranking if it uses rank normalization or its variants to process data.
\end{definition}
The two mainstream methods for optimization benchmarking, non-parametric hypothesis testing and Bayesian inference, are both considered to utilize relative ranking for data analysis.

\subsection{Non-parametric hypothesis testing and Bayesian inference}\label{subsec:non-parametric-hypothesis-testing-and-bayesian-inference}
Non-parametric hypothesis testing typically ranks algorithms based on the average of ranks they receive.
The primary objective of employing non-parametric testing techniques is to assess the statistical significance of the comparison results.
The general procedure, as outlined in~\cite{garcia_study_2009}, begins by testing whether there exists a significant performance difference among algorithms, typically utilizing Friedman's test~\cite{pereira2015overview} or Iman and Davenport's test~\cite{iman1980approximations}.
If such a difference is established, a reference control algorithm, like the algorithm with the best average rank, is selected.
Subsequently, Bonferroni-Dunn's procedure is applied to calculate the Critical Difference (CD) value.
This CD value is then used to create a visual representation indicating which algorithms significantly differ from the reference algorithm.
Additionally, the analysis is supplemented with a table displaying multiple adjusted $p$-values from algorithm comparisons, providing valuable reference information for the reader.

However, several weaknesses of Non-parametric hypothesis testing are pointed out~\cite{benavoli2017time, iman1980approximations}, especially regarding the misinterpretations of $p$-values~\cite{ioannidis2005most,wasserstein2016asa,goodman2008dirty,greenland2016statistical}.
Subsequently, Bayesian inference became an alternative tool.

The basic idea of Bayesian inference in algorithm benchmarking, as presented in~\cite{rojas-delgado_bayesian_2022}, is as follows.
First, an a probability model is chosen, which has several parameters to be determined, determining the probability of various outcomes occurring in pairwise comparisons of algorithms.
For example, the Bradley-Terry model~\cite{bradley1952rank} has $n$ parameters $(\theta_1, \theta_2, \dots, \theta_n)$ corresponding to $n$ algorithms respectively, and then the probability of $a_i\succ a_j$ is calculated by $P(a_i\succ a_j)=\frac{\theta_i}{\theta_i+\theta_j}$.
Then, a initial parameter vector is chosen as the prior distribution, and a new parameter vector is computed by maximizing the posterior probability using the results of the pairwise comparisons in the dataset.
Finally, substitute the obtained parameters into the probability model to derive the probability of assessing the algorithm's quality.

In the next section, we will demonstrate through a case study that these two methods may lead to NIIA issues.

    \section{The Issue of Relative Ranking}\label{sec:the-issue-of-relative-ranking}
    This section introduces the core issue addressed in this paper, which is the NIIA problem caused by the adoption of relative ranking.

\subsection{The NIIA issue in algorithm benchmarking}\label{subsec:the-niia-issue-in-algorithm-benchmarking}
\begin{figure}[!t]
    \centering
    \includegraphics[width=3.5in]{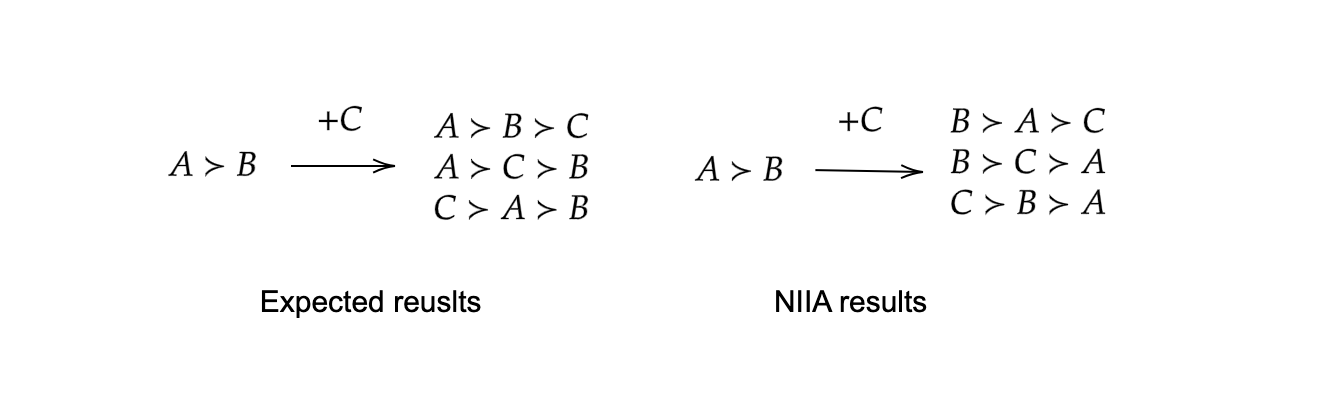}
    \caption{A schematic diagram illustrating the NIIA issue.
    The irrelevant algorithm $C$ unexpectedly influenced the comparison result of $A$ and $B$.}
    \label{fig:iia}
\end{figure}
Fig.~\ref{fig:iia} illustrates the NIIA issue in algorithm benchmarking.
When we compare algorithms $A$ and $B$, and obtain a result, let's assume $A\succ B$.
Then, when we introduce algorithm $C$ into the algorithm set, and provided that the analysis method is accurate, we would expect one of the following outcomes: $C\succ A\succ B$, $A\succ C\succ B$, or $A\succ B\succ C$.
However, if the results show that the addition of $C$ results in $B \succ A$, we can say that the analysis method exhibits the NIIA phenomenon.
Below, we provide a formal symbolic definition of the NIIA phenomenon.

Let the dataset $M$ represent the performance metrics of the algorithm set $\mathbb{A}$ on the benchmark function set $\mathbb{F}$.
Set $\mathbb{E}$ is a subset of $\mathbb{A}$, and $M(\mathbb{E})$ represents the projection of $M$,
with respect to the set $\mathbb{E}$.
Analysis method $T$ takes a dataset as input and generates some conclusions.
$T_{A,B}$ represents the comparison results of methods regarding algorithms $A$ and $B$.
Then $T$ is considered to exhibit the NIIA issue if and only if there exists $A,B,S$, $\{A, B\}\subseteq \mathbb{E} \subset \mathbb{A}$, such that,
\begin{equation}
    T_{A,B}(M) \neq T_{A,B}(M(\mathbb{E})).
\end{equation}

\subsection{A case of NIIA}\label{subsec:a-case-of-niia}
Provided that there are 500 problems and 100 algorithms.
In this context, we designate two specific algorithms as $A$ and $B$, while the remaining algorithms are labeled as $C_1$ through $C_{98}$.
For the first 100 problems, we assume the following performance metrics: Algorithm $A$ achieves a score of 1, Algorithm $B$ obtains a score of 99, Algorithm $C_1$ receives a score of 100, and for any $C_i$ (where $i$ ranges from 2 to 98), its score is equal to $i$.
In contrast, for the subsequent 400 problems, we assume the following performance metrics: Algorithm $A$ is assigned a score of 99, Algorithm $B$ is attributed a score of 98, Algorithm $C_1$ maintains a score of 100, and for any $C_i$ (where $i$ ranges from 2 to 98), its score becomes $i-1$.
In this scenario, we operate under the assumption that we can accurately collect and record these performance metrics for each experimental case, and that the lower scores indicate better performance.

To illustrate the NIIA issue, we construct two datasets using the example above.
The first dataset only observes comparisons between $A$, $B$, and $C_1$.
Its $M$ matrix is as follows:
\begin{equation}
    \bordermatrix{
        & I_1 & I_2 & \cdots & I_{100} & I_{101} & I_{102} & \cdots & I_{500}  \cr
        A & 1 & 1 & \cdots & 1 & 99 & 99 & \cdots & 99 \cr
        B & 99 & 99 & \cdots & 99 & 98 & 98 & \cdots & 98\cr
        C_1 & 100 & 100 & \cdots & 100 & 100 & 100 & \cdots & 100
    }.
\end{equation}

The second dataset, built upon the first dataset, includes all algorithms:
\begin{equation}
    \bordermatrix{
        & I_1 & I_2 & \cdots & I_{100} & I_{101} & I_{102} & \cdots & I_{500}  \cr
        A & 1 & 1 & \cdots & 1 & 99 & 99 & \cdots & 99 \cr
        B & 99 & 99 & \cdots & 99 & 98 & 98 & \cdots & 98\cr
        C_1 & 100 & 100 & \cdots & 100 & 100 & 100 & \cdots & 100 \cr
        C_2 & 2 & 2 & \cdots & 2 & 1 & 1 & \cdots & 1 \cr
        C_3 & 3 & 3 & \cdots & 3 & 2 & 2 & \cdots & 2 \cr
        \vdots & \vdots & \vdots & \vdots & \vdots & \vdots &\vdots & \vdots & \vdots \cr
        C_{97} & 98 & 98 & \cdots & 98 & 97 & 97 & \cdots & 97 \cr
    }.
\end{equation}

\begin{table}
    \begin{center}
        \caption{Result of non-parametric hypothesis testing.}
        \label{tab:npht}
        \begin{tabular}{c|cccccc}
            \hline
            & \makecell{Friedman\\ p-value} &
            \makecell{critical\\ difference\\($\alpha=0.001$)} & \makecell{average\\ rank\\ of A}
            & \makecell{average\\ rank\\ of B} & result\\
            \hline
            Dataset 1 & 3.94e-183 & 0.227 & 1.8 & 1.2 & $B\succ A$ \\
            Dataset 2 & 0.0 & 8.105 & 79.4 & 98.2 & $A \succ B$ \\
            \hline
        \end{tabular}
    \end{center}
\end{table}

First, we analyze these two datasets separately using non-parametric hypothesis testing method.
Friedman's test is adopted to assess overall differences, and the Bonferroni-Dunn's procedure is used to calculate the critical difference (CD) values.
The results are presented in Table~\ref{tab:npht}.
For dataset 1, the $p$-value obtained from Friedman's test falls below the significance level of $\alpha=0.001$, signifying a significant difference in algorithm performance.
The difference in the average rankings between Algorithm $A$ and Algorithm $B$ exceeds their respective CD value, which establishes the statistical significance of $A\succ B$.
However, when the same data and analytical process were applied to dataset 2, a contrary conclusion was reached, which is also deemed statistically significant.
From this results, it can be seen that analysis based on the average rank may lead to the NIIA issue, even when these results are considered statistically significant under non-parametric hypothesis testing.

Then, we explore the scenario of using Bayesian inference.
The experiment adopts the Bradley-Terry model, uses Zermelo's algorithm~\cite{zermelo1929berechnung} along with Newman's improved formula~\cite{newman2023efficient} to calculate the parameters.
The results are presented in Table~\ref{tab:bayes}, which reveals that Bayesian inference also encounters the NIIA issue.
\begin{table}
    \begin{center}
        \caption{Result of Bayesian inference.}
        \label{tab:bayes}
        \begin{tabular}{c|cccccc}
            \hline
            & $\theta_A$ & $\theta_B$ & $P(A\succ B)$ & $P(B\succ A)$ & result\\
            \hline
            Dataset 1 & 0.5 & 2.0 & 0.2 & 0.8 & $B\succ A$ \\
            Dataset 2 & 1.31e-4 & 1.01e-7 & 0.9992 & 0.0008 & $A\succ B$ \\
            \hline
        \end{tabular}
    \end{center}
\end{table}

The results indicate the potential presence of NIIA in the field of algorithm benchmarking.
In the following, we will analyze the causes behind it and propose potential solutions.

\subsection{The fundamental theorem about NIIA}\label{subsec:the-fundamental-theorem-about-niia}
The NIIA problem is not a phenomenon exclusive to algorithm benchmarking, it can occur in many scenarios that use relative orders, such as in the field of voting theory.
The main research focus of voting theory is aggregating individual rankings of candidates from each voter into a collective ranking that represents the overall consensus.
If we consider algorithms as candidates and benchmark functions as voters, the input of voting theory can be thought of as a data structure equivalent to the matrix $R$, where each column of data represents a permutation.

Equating voting theory and algorithm benchmarking in the manner described above can provide us with the following three insights:
\begin{enumerate}
    \item {Since the mathematical problems abstracted from the issues they each deal with are equivalent, most of their solutions can be interchangeable.
    For example, we can apply voting methods like Schulze~\cite{schulze2018schulze}, Ranked pairs~\cite{tideman1987independence}, Kemeny-Young~\cite{kemeny_mathematics_1959,young_consistent_1978}, Borda count~\cite{emerson_original_2013} to algorithm benchmarking, or use non-parametric hypothesis testing, Bayesian inference, deep statistics~\cite{eftimov_novel_2017} to rank candidates.
    In fact, non-parametric hypothesis testing and Borda's method share the same ranking strategy, the former essentially investigates the confidence level based on the conclusions of the latter.}
    \item {They face similar problems.
    For instance, the NIIA problem introduced in this paper also exists in voting theory.
    Liu~\cite{yan2022paradox} pointed out the existence of the cyclical ranking paradox in algorithm benchmarking, which is known as the Condorcet paradox in voting theory.}
    \item {They are constrained by the same underlying theorems, such as famous Arrow's impossibility theorem~\cite{arrow1951social}.}
\end{enumerate}

Arrow's impossibility theorem asserts that no voting system can simultaneously satisfy the criteria of Universal Domain, Non-Dictatorship, Pareto Principle, and Independence of Irrelevant Alternatives when there are three or more alternatives.
The accurate description and rigorous mathematical proof of Arrow's impossibility theorem can be found in~\cite{arrow1951social}.
For algorithm benchmarking, this theorem signifies that achieving a perfect solution to the NIIA issue under the premise of using relative rankings is impossible.
Fortunately, employing relative rankings is not a mandatory choice in algorithm benchmarking, and avoiding their use in this context is feasible because we have access to quantifiable raw data.
Below, we provide a straightforward explanation of how relative rankings can lead to issues.

\subsection{The Cause of NIIA}\label{subsec:the-cause-of-niia}
\begin{figure}[!t]
    \centering
    \includegraphics[width=3.5in]{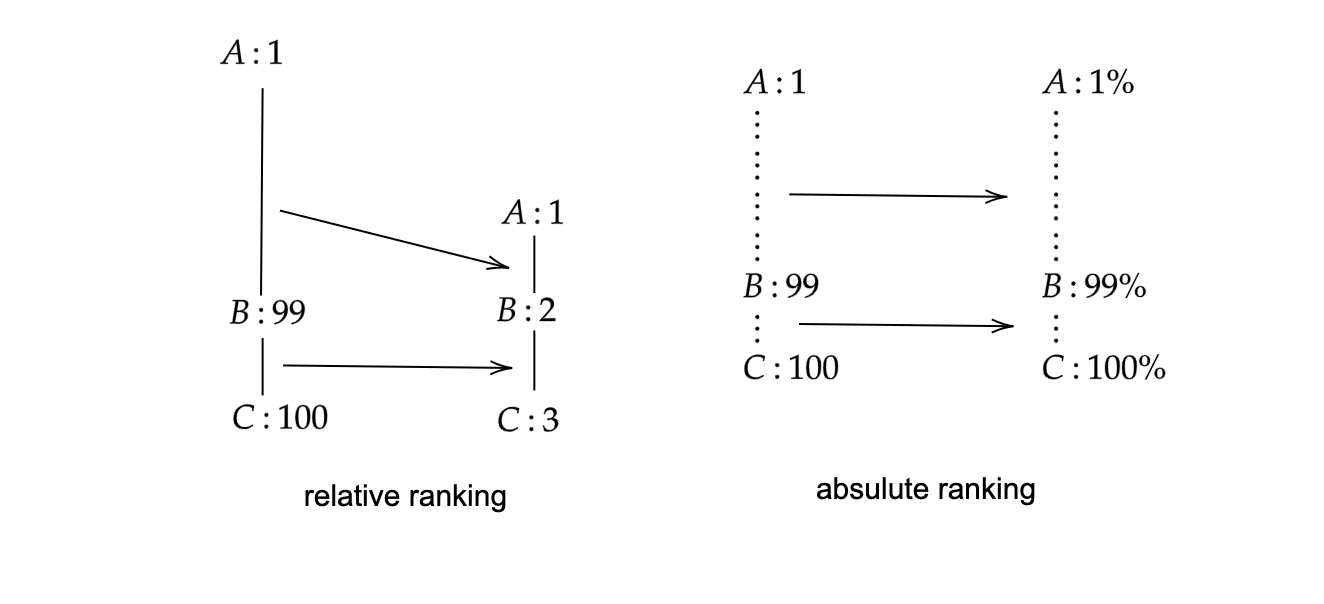}
    \caption{The erasure effect generated by using relative rankings.
    Algorithms positioned adjacently, regardless of the magnitude of their differences in raw data, have their gaps set to 1 after rank normalization.
    On the contrary, absolute ranking allows the performance differences between adjacent algorithms to become evident in their scores after normalization.}
    \label{fig:erase}
\end{figure}

By comparing the original data with the data after rank normalization, we observed that rank normalization did not induce any changes for Dataset 2.
In contrast, for Dataset 1, especially in the first 100 columns, rank normalization transformed the original values (1, 99, 100) to (1, 2, 3), which changes a lot.
This led to a misalignment between the data characteristics of $R$ and $M$, the difference between $A$ and $B$ is much greater than the difference between $B$ and $C_1$, while in $R$, the difference between $A$ and $B$ is the same as the difference between $B$ and $C_1$.
This suggests that during the process of rank normalization, detailed information about differences between adjacent rankings is erased.
Fig.~\ref{fig:erase} provides a visualization of this process.

Many data processing procedures exhibit similar phenomenon of information loss.
To describe it conveniently, we propose a definition as,

\begin{definition}[erasure effect]
    The erasure effect refers to a phenomenon in which data cannot be restored to its original value after processing.
    Accordingly, the data processing that produces the erasure effect is called an erasure process.
\end{definition}

With this definition, we can describe the phenomenon shown in Fig.~\ref{fig:erase} as the erasure effect of rank normalization.
It is important to note that not all erasure effects are harmful, for example, under the assumption that the data is generated with Gaussian error, finding the mean value of the data can be regarded as erasing the Gaussian error, which is a beneficial erasure process.
Normalization is also a beneficial erasure process as it erases the differences in the numerical scales of algorithm performance across different problems.

However, as we can see from Fig.~\ref{fig:erase}, the reason for the different conclusions drawn from dataset 1 compared to dataset 2 is that the erasure effect of rank normalization in Dataset 1 is stronger, to the extent that the rank normalized data $R$ differs significantly from the actual data $M$.

\subsection{Main idea of tackling NIIA}\label{subsec:main-idea-of-tackling-niia}
\begin{figure}[!t]
    \centering
    \includegraphics[width=3.5in]{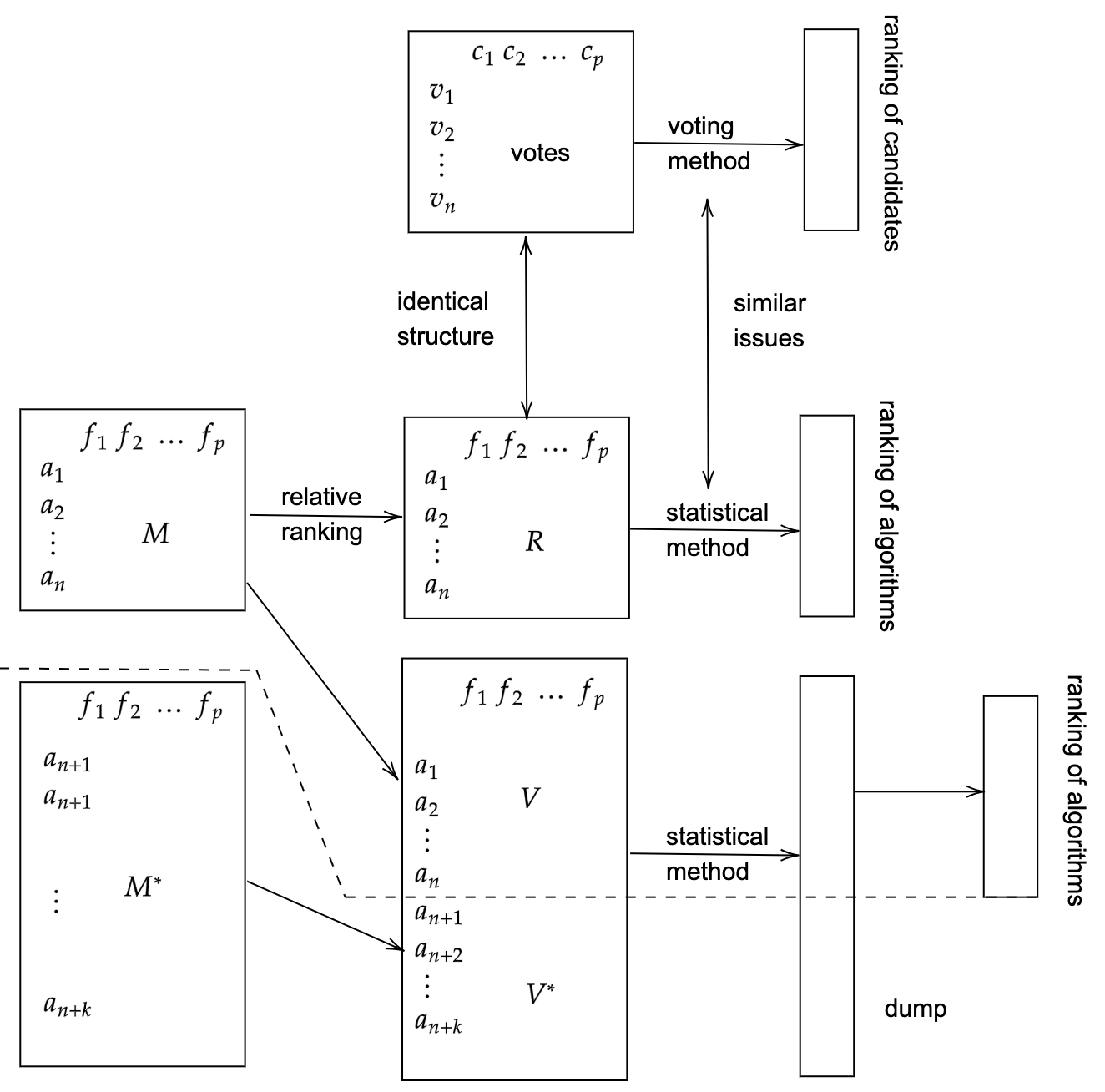}
    \caption{The main idea of tackling NIIA, where $M$, $M^*$, $R$, $V$, $V^*$ are all matrices.}
    \label{fig:idea}
\end{figure}
The results for dataset 2 were found to be more plausible than those for dataset 1 because it produced fewer erasure effects on the difference in algorithm performance.
Comparing the two datasets, we find that dataset 2 adds a large number of algorithms to dataset 1, which show a certain homogeneity in performance.
Since the numerical difference between $A$ and $B$ is large, there is a higher probability that the additional algorithms lie between $A$ and $B$, resulting in a much lower ranking of $A$ than $B$.

Inspired by this phenomenon, we have developed a method to reduce the excessive erasure effect of rank normalization.
Specifically, by introducing additional auxiliary algorithms like $C_2$ through $C_{98}$, the scale of differences between adjacent ranking algorithms can be reflected in their rank values.
A schematic representation of this approach is shown in Fig.~\ref{fig:idea}, which also demonstrates the connection between the algorithm benchmarking and voting theory.

A reasonable conjecture is that the more auxiliary algorithms are introduced, the less the ranking of algorithms is influenced by the other algorithms in the original algorithm set $R$, leading to more accurate results.
Let $\mathbb{S}$ be a set containing all possible algorithms, we can derive a new metric for evaluating algorithm performance as,
\begin{equation}
    \label{eq:va}
    v(a) = \frac{\pi_\mathbb{S}(a)}{|\mathbb{S}|},
\end{equation}
where $\pi_\mathbb{S}(a)$ is the rank of $a$ in $\mathbb{S}$ ordered by its performance metric and $\mathbb{S}$ is the cardinality of $\mathbb{S}$.

When adding or removing $k$ algorithms from the original algorithm set $\mathbb{A}$, the change of $\pi_\mathbb{S}(a)$ is limited to a maximum of $k$.
Given that the denominator approaches positive infinity, we can thus say that the new rank $v(a)$ can circumvent the NIIA problem.
Since this rank is not a relative rank in $\mathbb{A}$ but an absolute value, we shall name it the ``absolute rank''.

\begin{figure}[!t]
    \centering
    \includegraphics[width=3.5in]{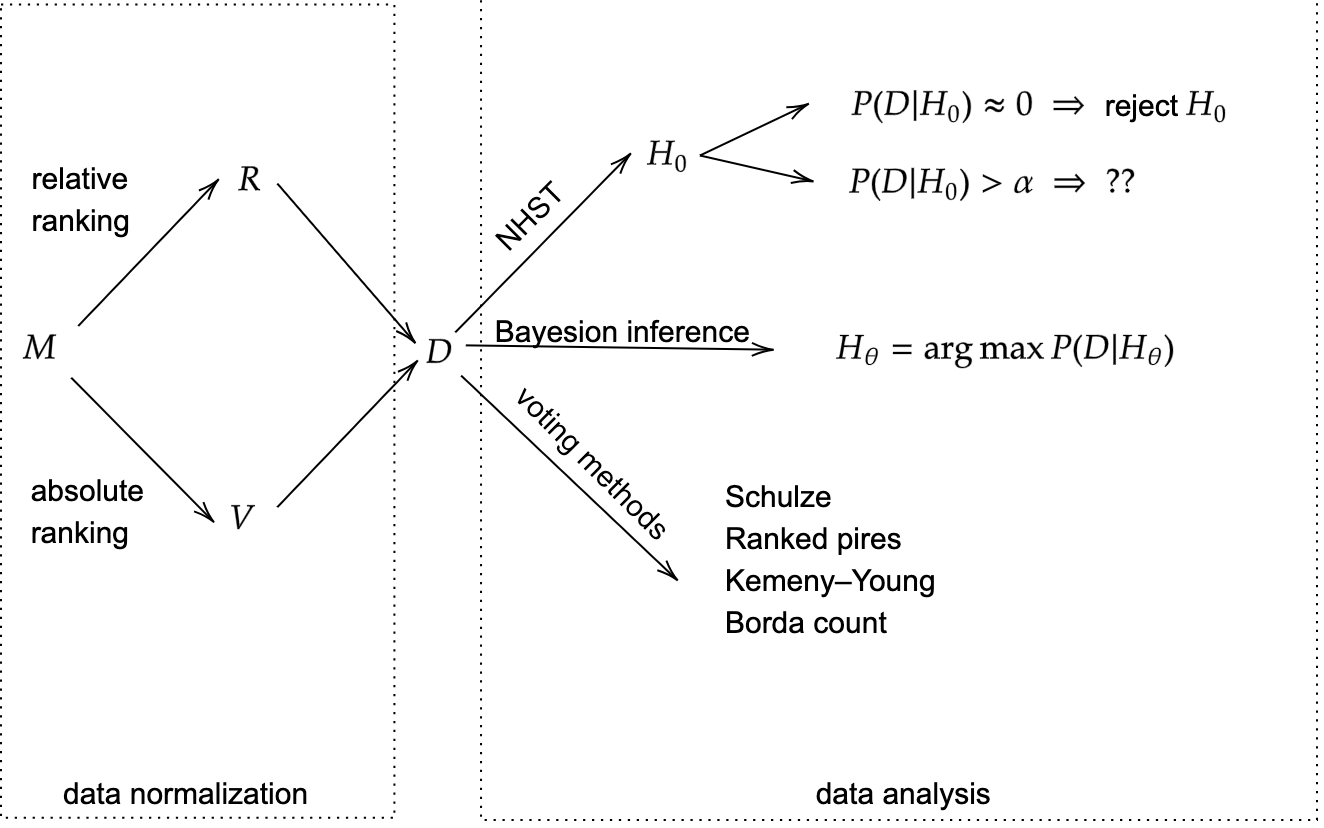}
    \caption{Illustration of the relationship between absolute ranking and other concepts.
    Absolute ranking is a process used to replace relative rankings for normalization, and its results are passed to data analysis methods for evaluation of optimization algorithms.}
    \label{fig:relation}
\end{figure}

To better illustrate the distinction between the ``absolute ranking'' and other conventional methods, we have depicted Fig.~\ref{fig:relation}.
It shows that the ``absolute ranking'' is not a data analysis procedure but rather an normalization before that, which is a better alternative to replace the original relative ranking in the context of algorithm benchmarking.

The mathematical model, computational methods, and applications of the ``absolute ranking'' will be discussed in the following sections.

    \section{Absolute Ranking}\label{sec:absolute-ranking}
    This section introduces the modeling of absolute ranking and discusses its computational methods.

\subsection{Modeling of absolute ranking}\label{subsec:modeling-of-absolute-ranking}
From the Section~\ref{sec:preliminaries}, we learn that the performance of algorithm $a$ on the benchmarking function $f$ is determined by its search path $X_{c*r}=X(a)$.
Provided that the domain of the benchmark function $f$ is $\Omega$, then the set containing all possible $X$ is:
\begin{equation}
    \mathbb{S} = \Omega^{c*r} = \{X_{c*r}|X_{ij}\in\Omega, i\in\{1,2,\dots,c\}, j\in\{1, 2, \dots, r\}\}
\end{equation}
Clearly, any search sequence obtained by an algorithm on the benchmark function $f$ belongs to the set $\Omega^{c*r}$, and the experimental results fall within $f(\Omega^{c*r})$, with their metrics belonging to $m(\Omega^{c*r})$.

If we replace the algorithm $a$ with its search sequence $X(a)$, Equation~\ref{eq:va} can be represented as:
\begin{equation}
    \label{eq:vx}
    v_f(X) = \frac{\pi_{m(\Omega^{c*r})}(m(X))}{|\Omega^{c*r}|} = \mathop{P}\limits_{x\sim U(\Omega^{c*r})}(m(x)\leq m(X))
\end{equation}
The probabilistic form is used to avoid difficulties in understanding caused by infinite denominators.
It is an equivalent but more easily comprehensible representation of the original expression.

Then the function to compute the absolute rank of a metric can be represented as,
\begin{equation}
    \label{eq:vt}
    v_f(t) = \mathop{P}\limits_{x\sim U(\Omega^{c*r})}(m(x)\leq t),
\end{equation}
where the input $t$ is the result of $m(a)$, and $v(t)$ is the absolute rank of $a$.

Using this function, we can replace the matrix $R$ composed of relative ranks with the one composed of the absolute ranks, represented by $V$, where
\begin{equation}
    \label{vij}
    V_{ij} = v_{f_j}(M_{ij}).
\end{equation}
Based on the previous analysis, replacing $R$ with $V$ will eliminate the NIIA issue.

The key step in the entire process is how to obtain an absolute ranking function like Equation~\ref{eq:vt} through a benchmark function.
This process is defined as follows,
\begin{equation}
    T_m(f) = v_f = (t)\to \mathop{P}\limits_{x\sim U(\Omega^{c*r})}(m(x)\leq t)
\end{equation}

To demonstrate the calculation process of $T_m$, we have employed a simple and commonly used metric, $m_0$, with the following definition,
\begin{equation}
    m_0(X) = \frac{1}{r}\sum_{j=1}^r \mathop{\min}\limits_{i=1}^c\{f(X_{ij})\}
\end{equation}

Let $\mathcal{X}$ be a random variable following a uniform distribution over the set $\Omega^{c*r}$, then $\mathcal{M}=m_0(\mathcal{X})$ is also a random variable.
Equation~\ref{eq:vt} is essentially calculating the cumulative distribution function (CDF) of the random variable $\mathcal{M}$.

When $c=1$ and $r=1$, Equation~\ref{eq:vt} can be written as,
\begin{equation}
    \label{eq:v1t}
    v_1(t) = \frac{\mathcal{L}^d \{\Omega(t-)\}}{\mathcal{L}^d \{\Omega\}},
\end{equation}
where $d$ is the dimension of $\Omega$, $\mathcal{L}^d$ is Lebesgue measure and $\Omega(t-) = \{x\in\Omega|f(x)\leq t)\}$.
Let $f_{min} = \min\{f(\Omega)\}$ and $f_{max}=\max\{f(\Omega)\}$, then we have $v(f_{min})=0$ and $v(f_{max})=1$.
Equation~\ref{eq:v1t} is essentially the CDF of the random variable $\mathcal{Y}=f(\mathcal{X})$.

When $r=1$ and $c>1$, Equation~\ref{eq:vt} calculates the CDF of the minimum value of multiple independent and identically distributed random variables $\mathcal{Y}_{min} = \min\{\mathcal{Y}_1, \mathcal{Y}_2, \dots, \mathcal{Y}_c\}$,
\begin{equation}
    \label{eq:vct}
    v_c(t) = 1-(1-v_1(t))^c.
\end{equation}

Furthermore, when $r>1$, the probability density of $\mathcal{M}$ can be obtained through $(r-1)$ convolutions.
As $r$ approaches infinity, the distribution of $M$ approaches a normal distribution with the mean and variance being the mean and variance of $\mathcal{Y}_{min}$, respectively.

Up to this point, we have transformed the calculation of $T_m$ into the determination of $v_1$.
In the following two subsections, we will introduce both analytical and sampling methods for solving $v_1$.

\subsection{Analysis of the absolute ranking}\label{subsec:analysis-of-the-absolute-ranking}

Solving Equation~\ref{eq:v1t} through mathematical analysis is often challenging, especially in high-dimensional or complex function expression cases.
However, in certain special scenarios, mathematical analysis is feasible and, in such cases, it is preferred over sampling methods because it provides a more accurate and elegant solution.

Under low-dimensional and simple function conditions, Equation~\ref{eq:v1t} can be solved using Riemann integration.
In addition, there are two other methods for solving absolute ranking in specific cases.

\begin{figure}[!t]
    \centering
    \includegraphics[width=3.5in]{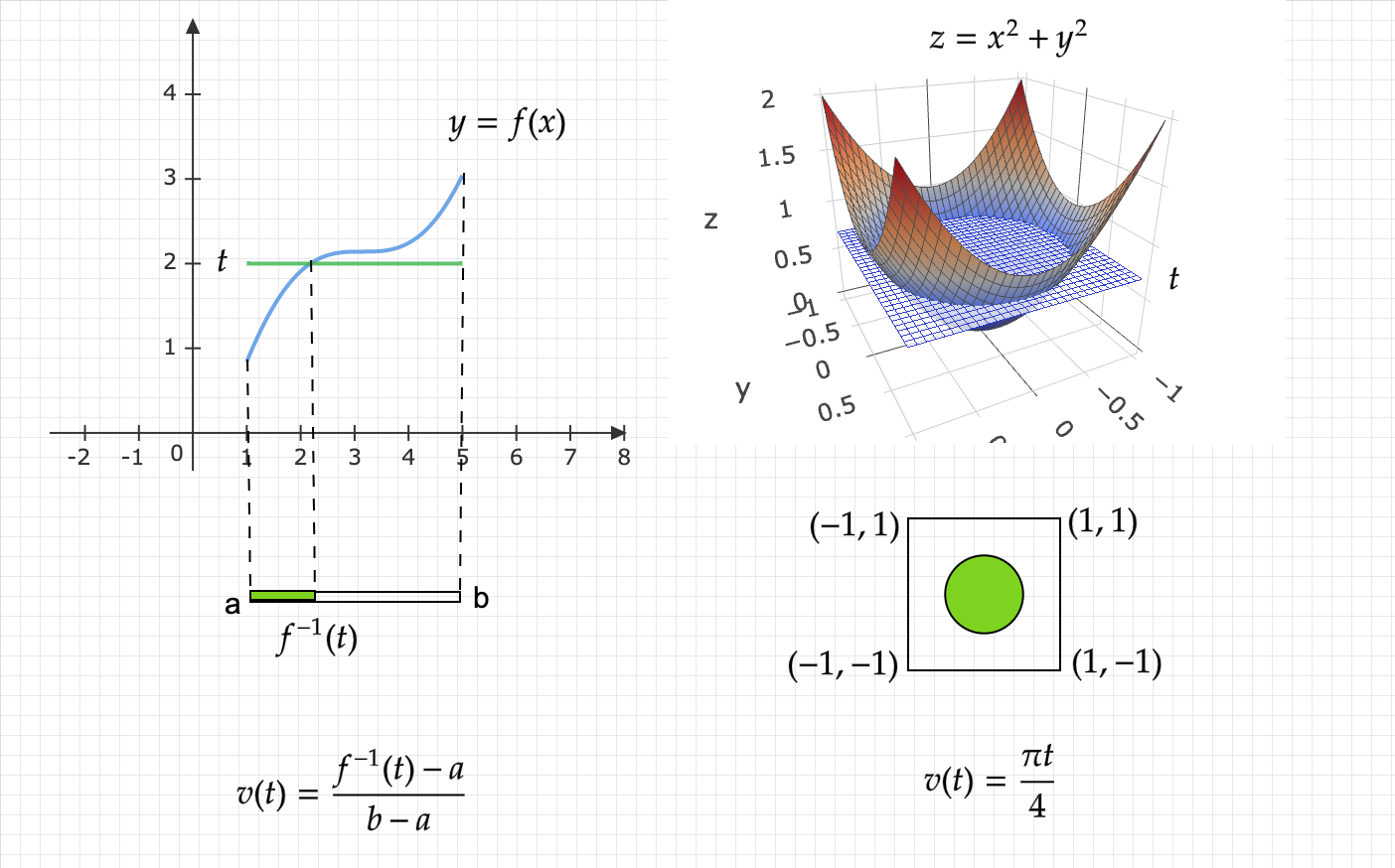}
    \caption{Illustration of calculating absolute rankings using Lebesgue measure.
    Absolute rank can be regarded as the proportion of the domain where the function values are less than a specific value.}
    \label{fig:vt}
\end{figure}
The first method is to directly utilize the Lebesgue measure formula for geometric shapes if the set $\Omega(t-)$ in space is a simple geometric shape.
Lebesgue measure can be interpreted as length, area, and volume in one-dimensional, two-dimensional, and three-dimensional spaces, respectively.
Fig.~\ref{fig:vt} illustrates two simple examples of using Lebesgue measure to calculate $v_1$.
In practice, the Sphere function is a common benchmark function, and its simplest form is:
\begin{equation}
    f(x)=\sum_{i=1}^d x_i^2,
\end{equation}
where $\mathcal{L}^d \{\Omega(t-)\}$ is a $d$-dimensional sphere of radius $\sqrt{t}$.
The Lebesgue measure formula for a $d$-dimensional sphere with a radius of $r$ is known as,
\begin{equation}
    \mathcal{L}^d \{ \text{sphere of radius } r \}=\frac{\pi^\frac{d}{2}}{\Gamma(\frac{d}{2}+1)}r^d.
\end{equation}
Provided that $\Omega=\{(x_1, x_2, \dots, x_d)|x_i\in[-w,w], i\in\{1,2,\dots,d\}\}$ and $\sqrt{t}<w$, we have:
\begin{equation}
    \label{eq:vtsp}
    v_1(t) = \frac{(\pi t)^\frac{d}{2}}{(2w)^d\Gamma(\frac{d}{2}+1)}.
\end{equation}

The second method involves making assumptions about the probability distribution.
Considering that Equation~\ref{eq:v1t} can be regarded as the CDF of $\mathcal{Y}$, we can obtain $v_1$ by assuming the distribution of $\mathcal{Y}$.

For example, when $\mathcal{Y}$ follows a uniform distribution on $[y_{min}, y_{max}]$, we have:
\begin{equation}
    \label{eq:vmms}
    v_1(t) = \int_{y_{min}}^t \frac{1}{y_{max}-y_{min}}dy = \frac{t-y_{min}}{y_{max}-y_{min}}.
\end{equation}
This result coincides exactly with the formula for representing Max-Min Scaling in Equation~\ref{eq:mms}.

Another common distribution is the Gaussian distribution, which assumes a probability density function of:
\begin{equation}
    q(y) = \frac{1}{\sqrt{2\pi}\sigma}e^{-\frac{(y-y_\mu)^2}{2\sigma^2}}.
\end{equation}
Then we have:
\begin{equation}
    v_1(t) = \Phi(\frac{y-y_\mu}{\sigma}),
\end{equation}
It includes a z-score normalization formula.

\begin{figure}[!t]
    \centering
    \includegraphics[width=3.5in]{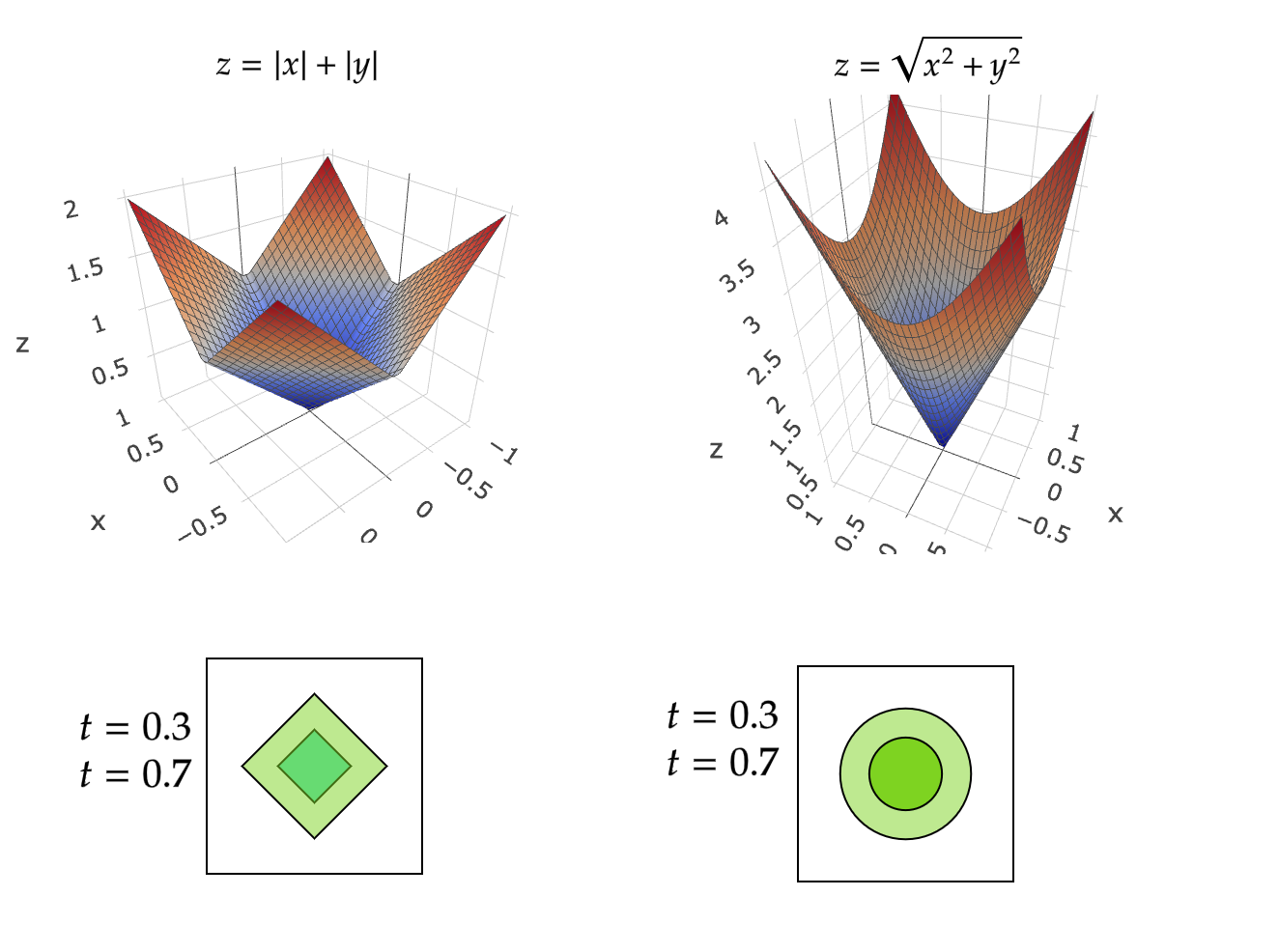}
    \caption{Similarity in contour plots of conic functions.
    In the figure, the light-shaded region represents the region in the domain where the function values are less than 0.7, while the dark-shaded area corresponds to 0.3.
    It can be observed that they exhibit similarity.}
    \label{fig:cone}
\end{figure}

While Max-Min Scaling and z-score normalization are commonly used tools in statistics, we have not found distinct functions that clearly adhere to their assumptions being used in optimization algorithm benchmarking.
However, through the analysis of certain common graphics, we have observed a regularity in the absolute ranking of a class of cone-shaped functions.
Specifically, as shown in Fig.~\ref{fig:cone}, when the function's graph is a pyramid or cone, let $S(t)$ be the geometric shape corresponding to $\Omega(t-)$, from Fig.~\ref{fig:cone} we can deduce that for any $t_1, t_2>0$, $S(t_1)$ is similar to $S(t_2)$ with a similarity ratio of $\lambda=\frac{t_1-y_{min}}{t_2-y_{min}}$, and consequently, their Lebesgue measures are in the ratio of $\lambda^d$.
Assuming there exists $t_0$ such that $\mathcal{L}^d(\Omega(t_0-))=1$, then we have:
\begin{equation}
    \label{eq:lcone}
    \mathcal{L}^d(\Omega(t-))=(\frac{t-y_{min}}{t_0-y_{min}})^d,
\end{equation}
provided $\Omega$ is $S(y_{max})$, substituting Equation~\ref{eq:lcone} into Equation~\ref{eq:v1t}, we get:
\begin{equation}
    \label{eq:vmmsd}
    v_1(t) = (\frac{t-y_{min}}{y_{max}-y_{min}})^d.
\end{equation}
Clearly, Equation~\ref{eq:vmms} is a special case of Equation~\ref{eq:vmmsd} when $d=1$.

Conic functions do have practical applications in algorithm benchmarking, such as function $F_5$ of the test functions for CEC 2005~\cite{suganthan2005problem}.

Up to this point, we have discussed the absolute ranking in the case of spherical and conic functions.
However, most functions are complex, so estimating the absolute rank through sampling is a more practical methodology.

\subsection{Estimate absolute ranking by sampling}\label{subsec:estimate-absolute-ranking-by-sampling}

In order to rationalize the sampling method, we introduce the following two assumptions first:
\begin{assumption}[Bounded finite discrete assumption(BFDA)]
    Limited by the level of human knowledge of the physical world, the precision of measurement or computational tools, there is always a precision threshold for real-world problems that makes research above this precision unattainable or meaningless, so that real-world problems can always be considered finite and discrete, and thus bounded.
\end{assumption}
\begin{assumption}[Local consistency assumption(LCA)]
    Based on the bounded finite dispersion assumption, the local consistency assumption assumes that the variation of the function's value in the minimum accuracy range is small enough to be ignored.
\end{assumption}
Based on these two assumptions, the fundamental idea of the sampling method to estimate $v_1$ is to randomly sample a large number of function values within the domain $\Omega$ and deduce the CDF of $\mathcal{Y}$ from these values.

The Monte Carlo method is typically an option for random search sampling.
However, in scenarios where the LCA holds, the Sobol method~\cite{sobol1967distribution} is considered superior to the Monte Carlo method.

\begin{figure}[!t]
    \centering
    \includegraphics[width=3.5in]{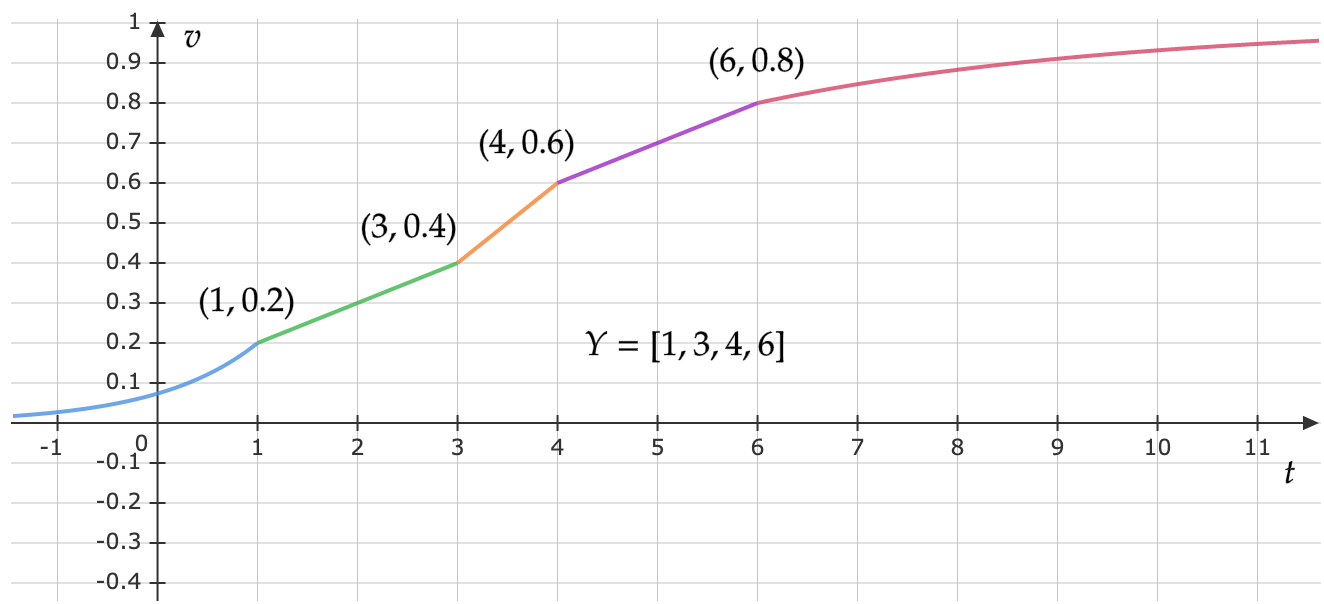}
    \caption{An example of a cumulative distribution function obtained from sampling results.}
    \label{fig:sobol}
\end{figure}

After $N=2^k$ samples by Sobol procedure, we get $N$ values of the benchmark function.
Denote the collection of these values as $Y$, and arrange its elements in ascending order as $y_1\leq y_2 \leq \dots \leq y_N$.
Then the absolute ranking function can be constructed as:
\begin{equation}
    \label{eq:vts}
    v_1(t)=\begin{cases}
               \frac{|\{y\in Y| y\leq t\}|}{N+1} & t\in Y\\
               v(y_i)+\frac{(t-y_i)(v(y_{i+1})-v_(y_i))}{y_{i+1}-y_i} & y_i < t < y_{i+1}\\
               \frac{e^t}{(N+1)e^{y_1}} & t<y_1\\
               1-e^{-\frac{\ln(n+1)}{y_n}t} & t>y_n
    \end{cases}.
\end{equation}

Formula~\ref{eq:vts} is a variant of the empirical distribution, ensuring that the entire function is strictly increasing and continuous.
The Fig.~\ref{fig:sobol} displays the graph of Equation~\ref{eq:vts} when $Y = [1,3,4,6]$, from which we can observe that the central part of the formula is processed using linear interpolation.
While this linear interpolation may appear less smooth with a small number of sample points, this drawback diminishes as the sample size becomes sufficiently large, and its advantage of being simple and computationally efficient becomes more pronounced.

Formula~\ref{eq:vts} assumes that the minimum and maximum values of the benchmark function within the domain are unknown.
Therefore, it employs exponential approximation to ensure that it can return a reasonable absolute ranking for any obtained value, whether it is very small or very large.
If the function's minimum value $y_{min}$ is known, $v(y)=0$ is hold for $y\leq y_{min}$, and $(y_{min}, 0)$ to $(y_1, \frac{1}{N+1})$ is connected with a straight line.
The same is done in a corresponding manner when the maximum value is known.

\begin{table}
    \begin{center}
        \caption{Sobol Sampling on Sphere function.}
        \label{tab:sobol}
        \begin{tabular}{ccc}
            \hline
            percentile & $y$ & $v(y)$\\
            \hline
            0 & 2.2748e-9 & 0.0002\%\\
            1st & 1.3206e-8 & 1.0001\%\\
            2nd & 1.5183e-8 & 2.0094\%\\
            3rd & 1.6506e-8 & 3.0510\%\\
            4th & 1.7525e-8 & 4.1165\%\\
            5th & 1.8371e-8 & 5.2111\%\\
            \hline
        \end{tabular}
    \end{center}
\end{table}

To verify the reliability of Sobol sampling and Formula~\ref{eq:vts}, we conduct Sobol sampling on a 10-dimensional spherical function.
Each dimension is confined to a sampling range from $\num{-e-4}$ to $\num{e-4}$, and set $N = 2^{20}$.
The results of the minimum value as well as the values within 1st to 5th percentiles are listed in TABLE~\ref{tab:sobol}.
These values are substituted into Equation~\ref{eq:vtsp} with d=10, and they are listed in the third column of TABLE~\ref{tab:sobol} in percentage form.

The results show that Sobol sampling and Equation~\ref{eq:vts}  can give a good approximation of the absolute ranking.

    \section{Applications}\label{sec:applications}
    Applying sampling methods in practice requires a trade-off between precision and cost.
When significant resources are required to compute the absolute ranking of a function, the responsibility for bearing this cost should be determined.
The perspective of this paper is that, since the absolute ranking depends solely on the objective function, the benchmark provider should be responsible for providing the absolute ranking function.

A benchmark function provider has a responsibility to offer information that helps users assess the meaning of the values obtained by their algorithms.
Typically, benchmark providers offer the reference minimum value of the benchmark function, but this is not enough.
For algorithms that can exactly retrieve this minimum value, the information is sufficient to prove their algorithm's superiority.
However, for those algorithms that do not find the optimal value, the difference between their computed value and the minimum value is affected by the numerical scale of the function, making the numerical significance unclear.
By providing a corresponding absolute ranking function for each benchmark function, users can gain a clearer understanding of the numerical meaning of their algorithm's results.

The ideal way to apply absolute ranking in practice is to establish a publicly accessible system that provides the implementation of common performance metrics, benchmark functions, and statistical methods as a public service to independent researchers.
This would allow researchers proposing new algorithms or algorithm improvements to focus their efforts on their research topics rather than diverting energy to collect algorithms from others to verify the effectiveness of their new algorithms.

Until such a public system is established, independent researchers may want to avoid or mitigate NIIA in their research at a lower cost.
In response to this demand, this paper offers some ideas for reference.

Firstly, areas closer to the optimal solution are considered more worthwhile to study than areas further away.
Define the range of the $i$-th dimension as $[x_i^* - \delta, x_i^* + delta]$, where $x^*$ is the optimal solution, thus obtaining a hypercube with a side length of $2\delta$ centered around the optimal solution.
Next, we need to select an appropriate $\delta$.
As shown in Fig.~\ref{fig:idea}, the addition of extra algorithms effectively mitigates the erasing effect of rank normalization.
However, not all algorithms have this effect.
When the added algorithm outperforms all algorithms in $\mathbb{F}$, the overall rank of all algorithms in $\mathbb{F}$ increases by 1 due to the addition of the new algorithm, which does not affect the algorithm's final rank.
When $\delta$ is too small, the performance of all algorithms in $\mathbb{F}$ is worse than the sampled data, and when $\delta$ is too large, the performance of all algorithms in $\mathbb{F}$ is more likely to be lower than the sampled data.
Therefore, we need to select an appropriate $\delta$ that covers the performance of sampling data as comprehensively as possible, better distinguishing these algorithms.

How to determine what is an ``appropriate'' search range?
An intuitive way is to observe whether the data is differentiating the algorithms in a dispersed manner.
This paper provides an empirical metric called the geometric mean difference, formulated as:
\begin{equation}
    \label{eq:e}
    e = \left( \mathop{\Pi}\limits_{i=2}^k (v_i-v_{i-1}) \right)^\frac{1}{k-1},
\end{equation}
where $v_1$ to $v_k$ are ordered distinct absolute ranks, from the smallest to the largest.

This gives us a two-step strategy.
First, with lower costs, roughly explore the region centered around the optimum point to find an appropriate search range close to the performance range of the algorithms in $\mathbb{F}$.
Then, with higher costs, conduct a detailed search within this range to obtain more accurate absolute orders for each algorithm.

\begin{table}
    \begin{center}
        \caption{Average error rate and average relative rank obtained on three functions in CEC2005}
        \label{tab:data}
        \begin{tabular}{c|ccc|c}
            \hline
            Algorithm & $f_{15}$ & $f_{16}$ & $f_{17}$ & ARR\\
            \hline
            BLX-GL50 & 400(10) & 93.49(2) & 109(2) & (4.67)\\
            BLX-MA & 269.6(7) & 101.6(6) & 127(8) & (7)\\
            CoEVO & 293.8(8) & 177.2(11) & 211.8(10) & (9.67)\\
            DE & 259(6) & 113(9) & 115(5) & (6.67)\\
            DMS-L-PSO & 4.854(1) & 94.76(3) & 110.1(3) & (2.33)\\
            EDA & 365(9) & 143.9(10) & 156.8(9) & (9.33)\\
            G-CMA-ES & 228(4) & 91.3(1) & 123(7) & (4)\\
            K-PCX & 510(11) & 95.9(4) & 97.3(1) & (5.33)\\
            L-CMA-ES & 211(3) & 105(7) & 549(11) & (7)\\
            L-SaDE & 32(2) & 101.2(5) & 114.1(4) & (3.67)\\
            SPC-PNX & 253.8(5) & 109.6(8) & 119(6) & (6.33)\\
            \hline
        \end{tabular}
    \end{center}
\end{table}

Below, we'll illustrate this process with a specific example.
The data in TABLE~\ref{tab:data} is derived from Table 13 in~\cite{garcia_study_2009}, with the addition of corresponding relative ranks and their corresponding Average Relative Ranks(ARR).
The scores calculated by Equation~\ref{eq:e} for each function are obtained and listed in TABLE~\ref{tab:area}, centered on the respective optimal point, at $2^{15}$ cost, with $\delta$ from 0.1 to 0.5, respectively.

\begin{table}
    \begin{center}
        \caption{Geometric mean difference gained by different range}
        \label{tab:area}
        \begin{tabular}{c|ccc}
            \hline
            $\delta$ & $f_{15}$ & $f_{16}$ & $f_{17}$ \\
            \hline
            0.1 & 0.02 & 0.00 & 0.00 \\
            0.2 & 1.70 & \textbf{3.44} & \textbf{4.57} \\
            0.3 & \textbf{2.29} & 0.74 & 0.55 \\
            0.4 & 0.88 & 0.57 & 0.07 \\
            0.5 & 0.65 & 0.73 & 0.03 \\
            \hline
        \end{tabular}
    \end{center}
\end{table}

From TABLE~\ref{tab:area}, we can see that the semi-side lengths of the optimal search regions for the three functions are 0.3, 0.2, and 0.2, respectively.
We sampled each function within their respective search area with a capacity of $2^{20}$, and the results are shown in TABLE~\ref{tab:result}.

\begin{table}
    \begin{center}
        \caption{Absolute rank obtained on three functions in CEC2005}
        \label{tab:result}
        \begin{tabular}{c|ccc|c}
            \hline
            Algorithm & $f_{15}$ & $f_{16}$ & $f_{17}$ & AAR\\
            \hline
            BLX-GL50 & 51.09 & 16.03 & 17.55 & 28.22\\
            BLX-MA & 9.89 & 23.43 & 32.83 & 22.05\\
            CoEVO & 14.98 & 97.32 & 93.42 & 68.54\\
            DE & 8.09 & 36.19 & 22.17 & 22.15\\
            DMS-L-PSO & 0.00 & 17.07 & 18.36 & 11.81\\
            EDA & 37.30 & 73.83 & 62.02 & 57.72\\
            G-CMA-ES & 4.20 & 14.32 & 29.08 & 15.87\\
            K-PCX & 89.80 & 18.04 & 10.27 & 25.70\\
            L-CMA-ES & 2.80 & 27.00 & 100.00 & 43.27\\
            L-SaDE & 0.00 & 23.03 & 21.44 & 14.82\\
            SPC-PNX & 7.32 & 32.16 & 25.51 & 21.66 \\
            \hline
        \end{tabular}
    \end{center}
\end{table}

The results in TABLE~\ref{tab:result} share some similarities with those in TABLE~\ref{tab:data}.
For example, the sets consisting of the top three, middle five, and bottom three algorithms are the same in both tables.
However, there are some differences in the rankings within those sets.

    \section{Conclusions}\label{sec:conclusions}
    The absolute ranking is developed to eliminate the NIIA issue which is an undesired phenomenon may occur in benchmarking optimization algorithms.
    This paper delves into the development of a mathematical model for absolute ranking, elucidates the process of computing absolute ranks through sampling, and offers practical recommendations and discussions regarding its real-world applications.

    The NIIA issue within the realm of algorithm benchmarking pertains to a scenario where the evaluation of algorithms is impacted by the inclusion of unrelated third-party algorithms.
    This phenomenon is often attributed to the usage of rank normalization techniques.
    Many existing methods employ rank normalization or similar data processing, which can potentially introduce unreliability due to the risk of NIIA.
    In this paper, an illustrative example is presented to demonstrate that two commonly used methods, non-parametric hypothesis testing and Bayesian inference, may exhibit NIIA problems when put into practice.

    To tackle this issue, the paper proposes a new ranking approach referred to as an absolute ranking, which carries more information about the objective function.
    Absolute ranking serves as an alternative to rank normalization and is designed not to be influenced by the performance of other algorithms, theoretically mitigating the occurrence of NIIA.
    This paper also provides analytical solutions for achieving the absolute ranking of two straightforward functions, namely the spherical function and the conical function.
    Additionally, it offers an estimation approach for applying the absolute ranking to any function through the utilization of sampling techniques.
    Finally, the paper delves into discussions on how these findings can be practically applied in real-world scenarios.

    There are several more in-depth topics that can be explored based on this work.
    Three main directions for future research are listed as below.
    \begin{enumerate}
        \item {Investigating the specific performance of different combinations of performance metrics and statistical methods under the absolute order scheme, as these factors can influence the entire algorithm benchmark testing process.}
        \item {Delving deeper into the analytical solutions of absolute orders for complex functions, including the study of absolute orders for basic elementary functions and the impact of various function compositions and transformations (such as translation and rotation) on absolute orders.}
        \item {Optimizing sampling methods, as Sobol sampling aims to sample evenly within a bounded domain, but in the context of optimization algorithms, regions with lower values are more valuable for sampling., which could lead to new, more targeted sampling methods.}
    \end{enumerate}

    Declaration of generative AI and AI-assisted technologies in the writing process

    During the preparation of this work the authors used Chat-GPT in order to
    improve language and readability.
    After using this tool, the authors reviewed and edited the content as needed and take full responsibility for the content of the publication.
    \bibliographystyle{IEEEtran}
    \bibliography{norm}

\end{document}